\documentclass[a4paper,12pt]{article}
\usepackage[psamsfonts]{amssymb}
\usepackage[psamsfonts]{eucal}

\usepackage{amsmath}
\usepackage{theorem}
\newtheorem{thm}{Theorem}
\newtheorem{prop}[thm]{Proposition}
\newtheorem{cor}[thm]{Corollary}
\newtheorem{lm}[thm]{Lemma}

\newtheorem{defin}{Definition}

\theorembodyfont{\rmfamily}

\theorembodyfont{\rmfamily}

\begin{document}

\title{Generalized pseudo-K\"ahler structures}

\author{Johann Davidov\thanks{Partially supported by "L.Karavelov" Civil Engineering Higher School, Sofia, Bulgaria
under contract No 10/2009}, Gueo Grantcharov, Oleg
Mushkarov\thanks{Partially supported by CNRS-BAS joint research
project {\it Invariant metrics and complex geometry, 2008-2009}},\\
Miroslav Yotov}

\date{}
\maketitle

\rm
\begin{abstract}
In this paper we consider pseudo-bihermitian structures -- pairs
of complex structures compatible with a
pseudo-Riemannian metric. We establish relations of these
structures with generalized (pseudo-) K\"ahler geometry and
holomorphic Poisson structures similar to that in the positive
definite case. We provide a list of compact complex surfaces which
could admit pseudo-bihermitian structures and give examples of
such structures on some of them. We also consider a na-turally
defined null plane distribution on a generalized pseudo-K\"ahler
$4$-manifold and show that under a mild restriction it determines
an Engel structure.

\end{abstract}

\section{Introduction}

Bihermitian structures have recently received a serious attention
due to their relations to supersymmetric sigma models in
theoretical physics and generalized geometry. However one of the
reasons they were introduced in \cite{AGG} was the observation
that the self-dual component of the Weyl tensor of an oriented
Riemannian 4-manifold determines a restriction on the number of
(local) complex structures compatible with the metric and the
orientation. The possibilities are 0, 1, 2, or $\infty$, if we
do not distinguish structures differing by sign. The bihermitian
structures thus arise naturally on 4-manifolds with 2 different
(up to sign) compatible complex structures. About 15 years earlier
than the paper \cite{AGG}, these structures appeared in the
physics literature  \cite{GHR}, where the target spaces of the
sigma-models with $(2,2)$-sypersymmetry were identified with
Riemannian manifolds admitting 2 compatible
complex structures satisfying additional differential
restrictions. An impulse for development of this topic in geometry
and string theory was the interpretation of bihermitian structures
in terms of the so-called generalized K\"ahler structures
\cite{Hit2, G1}, the latter being equivalent to the geometry
induced on the target of a $N=(2,2)$ supersymmetric sigma model
\cite{GHR, G}. This interpretation brought an important new
viewpoint for studying deformations of such structures and led to
a number of new examples \cite{G3, Goto}.

On a pseudo-Riemannian 4-manifold of neutral signature $(+,+,-,-)$
there are analogs for most of the notions in the Riemannian case.
In particular, compatible complex structures and self-duality are well defined, unlike the Lorentzian
case. Many results in the neutral setting are similar to results
in the Riemannian case but there are also important differences.

In this note we develop the notion of a pseudo-bihermitian structure which was considered also in
the physics literature \cite{L1}. We show that, in the same way as in the Riemannian case, it can
be related to (twisted) generalized pseudo-K\"ahler structures (Section 3) as well as to
holomorphic Poisson structures (Section 4). In Section 5 we show that the 3-dimensional complex
flag manifold $Fl$ carries a generalized K\"ahler structure. We also prove that any holomorphic
line bundle on $Fl$ is a holomorphic Poisson module with respect to a Poisson structure of a
special type. In Section 6 we provide a list of all compact complex surfaces which might carry
pseudo-bihermitian structures. It contains the list of bihermitian surfaces obtained in \cite{AGG}.
In Section 7 we adapt a construction of \cite{Hit1, G1} to find examples of pseudo-bihermitian
structures, which are collected in Proposition 10. Note that no Kodaira surface admits generalized
K\"ahler structures \cite{A, AGu}, but it admits a generalized pseudo-K\"ahler structure.

 We consider also some other differences between the Riemannian and the neutral setting. The
first one is related to the basic observation that on a 4-dimensional vector space two complex
structures $J_+$ and $J_-$ inducing the same orientation are compatible with a positive-definite
inner product iff $J_+J_-+J_-J_+=2pId$  for a constant $p$ with $|p|<1$. The same holds for
structures compatible with a split-signature inner product, but this time $|p|>1$. The difference
appears when the above identities are considered globally on a  4-manifold. If $p$ is a function
with $|p|<1$ at each point, then there always exists a unique conformal class of positive-definite
metrics compatible with $J_+$ and $J_-$. However we show in Section 7, Example 3 that there are
compact 4-manifolds admitting two such structures $J_+$ and $J_-$ with $|p|>1$ at every point which
are not compatible with a global pseudo-Riemannian metric, despite the fact that locally such a
metric always exists. Another difference comes from the fact that there is a naturally defined
null-plane distribution on any pseudo-bihermitian manifold, which is totally real with respect to
both complex structures. We show in Section 8 that, under a mild restriction this distribution, is
an Engel structure, which is a good analog of a contact structure in dimension four \cite{V}.

\medskip
\noindent {\bf Acknowledgements}: The authors express their gratitude to V.Apostolov for  helpful
discussions and comments on a preliminary version of this paper. Part of this work was done during
the visit of the first and the third-named authors at the Abdus Salam School of Mathematical
Sciences, GC University Lahore, Pakistan and the second named author's visit to the Institute of
Mathematics and Informatics at the Bulgarian Academy of Sciences. The authors thank the two
institutions for their hospitality.

\section{Pseudo-bihermitian structures}

In this section we consider the indefinite analog of bihermitian structures on 4-manifolds.
  An {\bf almost para-hypercomplex structure} on a smooth $4$-manifold
$M$ (also called an {\bf almost complex product} \cite{AS} or a {\bf neutral
almost hypercomplex structure} \cite{DW}) consists of three
endomorphisms $J_1,J_2,J_3$ of $TM$ satisfying the relations
\begin{equation}\label{3.2}
J_1^2=-J_2^2=-J_3^2=-Id, \hspace{.1in} J_1J_2=-J_2J_1=J_3
\end{equation}
of the imaginary units of the paraquaternionic algebra (split
quaternions). A metric $g$ on $M$ is called {\bf compatible} with
the structure $\{J_1,J_2,J_3\}$ if
\begin{equation}\label{3.3}
g(J_1X,J_1Y)=-g(J_2X,J_2Y)=-g(J_3X,J_3Y)=g(X,Y)
\end{equation}
(such a metric is necessarily of neutral signature $(+,+,-,-)$).
In this case we say that $\{g,J_1,J_2,J_3\}$ is an {\bf almost
para-hyperhermitian structure}. For any such a structure we define
three $2$-forms $\Omega_i$ setting $\Omega_i(X,Y) = g(J_iX,Y)$,
$i=1,2,3$. If the Nijenhuis tensors of $J_1,J_2,J_3$ vanish, the
structure $\{g,J_1,J_2,J_3\}$ is called {\bf para-hyperhermitian} and $(J_1, J_2, J_3)$ is called {\bf para-hypercomplex}. When
additionally the 2-forms $\Omega_i(X,Y) = g(J_iX,Y)$ are closed,
the para-hyperhermitian structure is called {\bf para-hyperk\"ahler}
(also called {\bf hypersymplectic} \cite{Hit0} and {\bf neutral hyperk\"ahler}
\cite{DW}).

Hypercomplex or para-hypercomplex structures can be obtained in
the following way. Consider a 4-manifold with two complex
structures $J_+$ and $J_-$ such
 that
\begin{equation}\label{anticom}
 J_+J_-+J_-J_+ = 2p Id
\end{equation}
for a function $p$.


Suppose that $|p|<1$ at each point. Then $J_+$,
$K=\displaystyle{\frac{1}{2\sqrt{1-p^2}}}[J_+, J_-]$,
$S=-\displaystyle{\frac{1}{\sqrt{1-p^2}} }(J_- +pJ_+)$ form an
almost hypercomplex structure (cf. e.g.\cite{L1}). Thus the
complex structures $J_+$ and $J_-$ are compatible with a positive
definite metric.

 If $|p|>1$ at every point, then $$J_+, ~K=\displaystyle{\frac{1}{2\sqrt{p^2-1}}}[J_+,
 J_-],~ S=-\displaystyle{\frac{1}{\sqrt{p^2-1}} }(J_- +pJ_+)$$
form an almost para-hypercomplex structure \cite{L1}. Hence by
\cite{DGMY1} there is a locally defined metric compatible with the
structure $\{J_+,K,S\}$. It is clear that the structure $J_-$ is
also compatible with this metric. Conversely, if the structures
$J_+$ and $J_-$ are compatible with a pseudo-Riemannian metric
$g$, so will be $K$ and $S$, hence $g$ is of neutral signature.
Note that, unlike the positive definite case, given $J_+$ and
$J_-$, such a metric may not exist globally (see Example 3 in
Section 7).

It follows from the above discussion that if $|p|\neq 1$ at every
point, then $J_+$ and $J_-$ yield the same orientation. This is a
consequence from the well-known fact that two non-collinear
(almost) complex structures on a $4$-manifold both compatible with
a pseudo-Riemannian  metric determine opposite orientations
exactly when they commute.

\begin{defin}If $J_+ \neq \pm J_-$ are complex structures on a $4$-manifold
compatible with a pseudo-Riemannian metric $g$ and if they yield
the same orientation, then $(g,J_+,J_-)$ is said to be a {\it
pseudo-bihermitian structure}. Such a structure is called {\it
strict} if $J_+\neq \pm J_-$ at every point.
\end{defin}

Note that if $(g,J_+,J_-)$ is a pseudo-bihermitian structure, then
$J_+$ and $J_-$ satisfy identity (\ref{anticom}) with
$p=-\frac{1}{2}g(J_+,J_-)$.

The following lemma is well-known in the positive definite case.
For the neutral case it is stated in \cite{L1} and proved in
\cite{L2} for generalized K\"ahler structures. For the sake of
completeness we provide a new proof, which works both in the
positive and neutral-signature cases.

\begin{lm}
Let $J_+$ and $J_-$ be complex structures on a $4$-manifold such
that $J_+J_-+J_-J_+ = 2p Id$ for $p=const$ and $|p|>1$. Then
$\{J_+,K,S\}$ is a para-hypercomplex structure.
\end{lm}
{\it Proof:} We have to prove that the almost product structures $K$ and $S$ are integrable. To do
this we shall use a local neutral metric $g$ compatible with the structure $\{J_+,K,S\}$. Then
$J_-$ is also compatible with $g$  and $p=-\frac{1}{2}g(J_+
 ,J_-)$. Denote by $F^{\pm}$ the K\"ahler $2$-form of $(g,J_{\pm})$.
  Then a standard formula for the Hermitian structure $(g,J_{\pm})$ gives:
\begin{equation}\label{nabla}
\begin{array}{c}
g((\nabla_{X}J_{\pm})(Y),Z)=(\nabla_XF^{\pm})(Y,Z)=\\[6pt]
\displaystyle{\frac{1}{2}}(dF^{\pm}(J_{\pm}X,Y,J_{\pm}Z)+dF^{\pm}(J_{\pm}X,J_{\pm}Y,Z)),
\end{array}
\end{equation}
where $\nabla$ is the Levi-Civita connection of $g$.

Since the dimension of the manifold is four, there is a unique
$1$-form $\theta_{\pm}$ (the Lee form) such that $dF^{\pm} =
\theta_{\pm}\wedge F^{\pm}$.
 Then
$$
\begin{array}{lll}
g((\nabla_{X}J_{\pm})(Y),Z)&=&g(X,Z)\theta_{\pm}(J_{\pm}Y)-g(J_{\pm}X,Z)\theta_{\pm}(Y)\\[6pt]
                           & &-g(X,Y)\theta_{\pm}(J_{\pm}Z)-g(J_{\pm}X,Y)\theta_{\pm}(Z)
\end{array}
$$
It follows that
$$
\begin{array}{lll}
2X(g(J_{+},J_{-}))&=&2g(\nabla_{X}J_{+},J_{-})+2g(J_{+},\nabla_{X}J_{-})=\\[6pt]
  & &-\theta_{+}([J_{+},J_{-}]X)+\theta_{-}([J_{+},J_{-}]X)
\end{array}
$$
Thus
\begin{equation}\label{dp}
  2d(g(J_+,J_-))=-(\theta_+-\theta_-)\circ[J_+,J_-]
\end{equation}
In view of the identity $2p=-g(J_{+},J_{-})$, the condition
$p=const$ leads to $\theta_+=\theta_-$ since
$[J_{+},J_{-}]=2\sqrt{p^2-1}K\neq 0$ at every point.  Then, using
the identity
  $S=-\displaystyle{\frac{1}{\sqrt{p^2-1}}}(J_- +pJ_+)$, we see that the fundamental $2$-form $F^S$ of
$S$ is a linear combination of  $F^{-}$  and $F^{+}$ with constant
coefficients. Hence $dF^{S} = \theta_+\wedge F^{S}$, so the Lee
form of $(g,S)$ is $\theta_+$. Let $F^K$ be the fundamental
$2$-form of $(g,K)$ and denote its Lee form by $\theta^K$.  Take a
$g$-orthogonal basis of tangent vectors $\{E_1,E_2,E_3,E_4\}$ with
$||E_1||^2=||E_2||^2=1$, $||E_3||^2=||E_4||^2=-1$. Set
$\varepsilon_i=||E_i||^2$, $i=1,2,3,4$. Then the identities
$dF^K=\theta^K\wedge F^K$ and
$$
\begin{array}{l}
\sum_{i=1}^4 \varepsilon_i dF^K(E_i,KE_i,Z)= 2\sum_{i=1}^4
\varepsilon_ig((\nabla_{E_i}K)(KE_i),Z).
\end{array}
$$
give
$$
\theta^K(Z)=-\sum_{i=1}^4\varepsilon_i[g((\nabla_{E_i}K)(E_i),KZ).
$$
for any tangent vector $Z$. Since $K=-J_+S$, we have
$$
\theta^K(Z)=-\sum_{i=1}^4\varepsilon_i[g((\nabla_{E_i}J_+)(SE_i),J_+SZ)-
\sum_{i=1}^4\varepsilon_i[g((\nabla_{E_i}S)(E_i),SZ).
$$
Using (\ref{nabla}) and the fact that $dF^{+}=\theta_+\wedge F^{+}$ one can easily see that the
first term on the right-hand side vanishes. The second term is $\theta^S(Z)$. Thus
$\theta^K=\theta^{S}=\theta_+$, therefore the structures $K$ and $S$ are integrable \cite{Kam}.
{\it q.e.d.}

\section{Generalized pseudo-K\"ahler structures}

Recall that a $H$-twisted generalized complex structure on a smooth
manifold $M$ is an endomorphism $I$ of the bundle $TM\oplus
T^{\ast}M$ satisfying the following conditions:

\noindent $(a)$ $I^2=-Id$,

\medskip

\noindent$(b)$ $I$ preserves the natural metric

\medskip

\noindent $<X+\xi,Y+\eta>=\frac{1}{2}(\xi(Y)+\eta(X)), \quad
X,Y\in TM,\quad \xi,\eta\in T^{\ast}M$

\medskip

\noindent$(c)$ the $+i$-eigensubbundle of $I$ in $(TM\oplus
T^{\ast}M)\otimes {\Bbb C}$ is involutive with respect to the
$H$-twisted Courant bracket defined by
$$ [X+\xi, Y+\eta]_H=
[X,Y]+L_X\eta-L_Y\xi-\frac{1}{2}(d\imath_X\eta-d\imath_Y\xi) +
\imath_Y\imath_X H,
$$
where $H$ is a closed $3$-form.

 The integrability condition $(c)$ is equivalent to vanishing of the Nijenhuis tensor
$$
N_H(A,B)=[A,B]_H-[IA,IB]_H+I[IA,B]_H+I[A,IB]_H,~ A,B\in TM\oplus
T^{\ast}M.
$$

The space of 2-forms $\Omega^2(M)$ acts on $TM\oplus T^*M$ as $e^b(X+\xi)=X+\xi+\imath_Xb$ for any
$b\in\Omega^2(M)$. Then the Courant bracket satisfies $[e^b(A),e^b(B)]_H=[A,B]_{H+db}$. In
particular if $I$ is a generalized complex structure, integrable with respect to the $H$-twisted
Courant bracket, then $J=e^{-b}Ie^{b}$ is a generalized complex structure, integrable with respect
to the $(H-db)$-twisted Courant bracket. So whenever $H$ is exact, $H=db$ for some $2$-from $b$,
the structure $I$ is called untwisted since the structure $J$ is integrable with respect to the
Courant bracket with  vanishing $3$-form.

Following M.Gualtieri \cite{G1, G} we introduce the following:
\begin{defin}
A (twisted) generalized pseudo-K\"ahler structure is a pair of commuting (twisted) generalized
complex structures $I_1, I_2 : TM\oplus T^*M\rightarrow TM\oplus T^*M$, such that the
$\pm1$-eigenspaces $L^{\pm}$ of $G=I_1I_2$ are transversal to $TM$ and the canonical inner product
on $TM\oplus T^*M$ is non-degenerate on $L^{\pm}$.
\end{defin}

Using the same proof as in \cite{G1, G}, we have

\begin{thm}\label{gK}
A $H$-twisted generalized pseudo-K\"ahler structure on a manifold $M$ is equivalent to a quadruple
$(g,J_+, J_-,b)$, where $g$ is a pseudo-Riemannian metric, $J_+$ and $J_-$ are $g$-Hermitian
complex structures, and $b$ is a 2-form such that $d^+F^+ = -d^-F^-=H + db$, where $F^{\pm}$ is the
K\"ahler form of $(g,J_{\pm})$ and  $d^{\pm}$ is the imaginary part of the
$\overline{\partial}$-operator of $J_{\pm}$.
\end{thm}

\section{Holomorphic Poisson structures}

In this section we prove an indefinite analog of the well-known
result \cite{Hit1} that a generalized K\"ahler manifold carries a
holomorphic Poisson structure. In fact, we have the following
slightly more general result.

\begin{thm}\label{Poisson}
Let $(M,g)$ be a pseudo-Riemannian manifold and let $J_+$, $J_-$
be two complex structures on $M$ compatible with $g$ and such that
$d^+F^+ = -d^-F^- $. Then $M$ admits a $J_+$-holomorphic Poisson
structure which vanishes iff $[J^+,J^-]=0$.
\end{thm}
{\it Proof:} Let $\Pi$ be the bivector field on $M$ determined by the endomorphism
$[J_+,J_-]-iJ_+[J_+,J_-]$ of $T^{\Bbb C}M$ and the complex bilinear extension of $g$. We shall
prove that $\Pi$ is a holomorphic Poisson field. To show that $\Pi$ is holomorphic we shall use the
Chern connection $D^+$ of the pseudo-Hermitian structure $(g,J_+)$.  It is defined by the identity
$g(D^{+}_XY,Z) = g(\nabla_XY,Z) -\frac{1}{2}dF^{+}(J_{\pm}X,Y,Z)$, where $\nabla$ is the
Levi-Civita connection of $g$. As in the positive case, $D^{+}$ is a Hermitian connection such that
the restriction of its $(0,1)$ part on the holomorphic tangent bundle is the
$\overline{\partial}$-operator of $J_{+}$.

In view of (\ref{nabla}) and the identity
$d^{\pm}F^{\pm}(X,Y,Z)=-dF^{\pm}(J_{\pm}X,J_{\pm}Y,J_{\pm}Z)$, we
have
$$
\begin{array}{lll}
2g((D^+_XJ_-)(Y),Z)=2g(D_X^+J_-Y,Z)+2g(D^+_XY,J_-Z)=\\[6pt]
2g((\nabla_XJ_-)(Y),Z)-dF^+(J_+X,J_-Y,Z)-dF^+(J_+X,Y,J_-Z)=\\[6pt]
dF^-(J_-X,Y,J_-Z)+dF^-(J_-X,J_-Y,Z)\\[6pt]
\hspace{3.5cm}-dF^+(J_+X,J_-Y,Z)-dF^+(J_+X,Y,J_-Z)=\\[6pt]
d^-F^-(X,J_-Y,Z)+d^-F^-(X,Y,J_-Z)\\[6pt]
\hspace{3.5cm}+d^+F^+(X,J_+J_-Y,J_+Z)+d^+F^+(X,J_+Y,J_+J_-Z)
\end{array}
$$

 Thus
\begin{equation}\label{Chern1}
\begin{array}{lll}
2g((D^+_XJ_-)(Y),Z)=-d^+F^+(X,J_-Y,Z)-d^+F^+(X,Y,J_-Z)\\[6pt]
\hspace{3.5cm}+d^+F^+(X,J_+J_-Y,J_+Z)+d^+F^+(X,J_+Y,J_+J_-Z)
\end{array}
\end{equation}
The $3$-form $d^+F^+$ has no $(3,0)$ and $(0,3)$-components, so
$$
\begin{array}{lll}
d^+F^+(A,B,C)=\\[6pt]
\hspace{1.5cm}d^+F^+(J_+A,J_+B,C)+d^+F^+(J_+A,B,J_+C)+d^+F^+(A,J_+B,J_+C)
\end{array}
$$
Applying this identity to the last two terms in (\ref{Chern1}) we
get

\begin{equation}\label{Chern2}
\begin{array}{lll}
2g((D^+_XJ_-)(Y),Z)=-d^+F^+(J_+X,J_-Y,J_+Z)-d^+F^+(J_+X,J_+J_-Y,Z)\\[6pt]
\hspace{3.5cm}-d^+F^+(J_+X,Y,J_+J_-Z)-d^+F^+(J_+X,J_+Y,J_-Z)
\end{array}
\end{equation}

Set $Q=[J_+,J_-]$. Then, since $D^+J_+=0$, we have

$$
\begin{array}{c}
2g((D^+_{X}Q)(Y),Z)-g((D^+_{J_+X}Q)(Y),J_+Z)=\\[6pt]
-g((D^+_{X}J_-)(Y),J_+Z)-g((D^+_{X}J_-)(J_+Y),Z)\\[6pt]
-g((D^+_{J_+X}J_-)(Y),Z)+g((D^+_{J_+X}J_-)(J_+Y),J_+Z)
\end{array}
$$
Applying (\ref{Chern1}) to the first and the second term, and
(\ref{Chern2}) to the third and the fourth term, we easily get

\begin{equation}\label{derP}
g((D^+_{X}Q)(Y),Z)-g((D^+_{J_+X}Q)(Y),J_+Z)=0
\end{equation}
As in \cite{AGG} and \cite{Hit1}, consider the form
$\Omega(X,Y)=g(QX,Y)$. The $(1,1)$-part of this form with respect
to $J_+$ vanishes since
$$
\begin{array}{r}
\Omega(J_+X,J_+Y)=-g(J_+^2J_-J_+X,Y)-g(J_-J_+^2X,J_+Y)=\\[6pt]
g(J_-J_+X,Y)+g(J_-X,J_+Y)=-\Omega(X,Y).
\end{array}
$$
Then the $(0,2)$-component of $\Omega$ is
\begin{equation}\label{omega}
\begin{array}{r}
\Omega^{(0,2)}(X,Y)=\Omega^{(0,2)}(X^{(0,1)},Y^{(0,1)})= \Omega(X^{(0,1)},Y^{(0,1)})=\\[6pt]

\frac{1}{4}[\Omega(X,Y)-\Omega(J_+X,J_+Y)]+\frac{1}{4}i[\Omega(J_+X,Y)+\Omega(X,J_+Y)]=\\[6pt]
\frac{1}{2}[\Omega(X,Y)+i\Omega(X,J_+Y)]=\frac{1}{2}[g(QX,Y)+ig(QX,J_+Y)]=\frac{1}{2}g(\Pi,X\wedge
Y)
\end{array}
\end{equation}
It follows that $\Pi$ is of type $(2,0)$ with respect to $J_+$.
Moreover, we have
$$
\begin{array}{lrr}
g(D^{+}_{X+iJ_+X}\Pi,Y\wedge Z)=
2(D^{+}_{X+iJ_+X}\Omega^{(0,2)})(Y,Z)=\\[6pt]
[g((D^{+}_{X}Q)(Y),Z)+ig((D^{+}_{X}Q)(Y),J_+Z)]\\[6pt]
\hspace{0.3cm}+i[g((D^{+}_{J_+X}Q)(Y),Z)+ig((D^{+}_{J_+X}Q)(Y),J_+Z)]=\\[6pt]
[g((D^{+}_{X}Q)(Y),Z)- g((D^{+}_{J_+X}Q)(Y),J_+Z)]\\[6pt]
\hspace{0.3cm}+i[g((D^{+}_{X}Q)(Y),J_+Z)+g((D^{+}_{J_+X}Q)(Y),Z)]
\end{array}
$$
Hence, by (\ref{derP}),  $ g(D^{+}_{X+iJ_+X}\Pi,Y\wedge
Z)=0~\mbox{for every}~X,Y,Z\in TM. $ This shows that
$D^{+}_{X+iJ_+X}\Pi=0$, therefore $\Pi$ is a holomorphic section
of the anti-canonical bundle $\Lambda^2T^{(1,0)}M$ of $(M,J_{+})$.

To prove that the Schouten-Nijenhuis  bracket $[\Pi,\Pi]$ vanishes, we note first that it is enough
to show that $[Re\Pi,Re\Pi]=0$. Indeed, since $\Pi$ is holomorphic, it is easy to see in local
holomorphic coordinates that $[\overline{\Pi},\Pi]=0$. Note also that $[Re\Pi,Im\Pi]=[Im\Pi,Re\Pi]$
since $Re\Pi$ and $Im\Pi$ are of degree $2$.  Thus, we have
$0=[\overline{\Pi},\Pi]=[Re\Pi,Re\Pi]+[Im\Pi,Im\Pi]$. Suppose that $[Re\Pi,Re\Pi]=0$. Then we get
$[Im\Pi,Im\Pi]=0$, hence $[\Pi,\Pi]=2i[Re\Pi,Im\Pi]$. Because $[\Pi,\Pi]$ is of type (3,0) and
purely imaginary, we conclude that  $[\Pi,\Pi]=0$

 According to (\ref{omega}), the endomorphism $Q$ of $TM$ corresponds to the bivector field $Re\Pi$ via the metric $g$.
Then, in view of \cite[Proposition 1.9]{Vais}, the equality $[Re\Pi,Re\Pi]=0$ is equivalent to
$$
\mathfrak{G} g((\nabla_{QX}Q)(Y),Z)=0,
$$
where $\mathfrak{G}$ means the cyclic sum over $X,Y,Z$ and
$\nabla$ is the Levi-Civita connection of $g$. To prove the latter
identity we use the fact that the Levi-Civita connection $\nabla$
and the Chern connection $D^{+}$ of $(g,J_+)$ are related by

$$
g(\nabla_XY,Z) = g(D^{+}_XY,Z)
-\frac{1}{2}d^{+}F^{+}(X,J_{+}Y,J_{+}Z).
$$
Set $P=J_+J_-+J_-J_+$. Then, by (\ref{Chern1}) we have
\begin{equation}\label{nablaQ}
\begin{array}{l}
2g((\nabla_XQ)(Y),Z)=2g((\nabla_XQY,Z)+2g(\nabla_XY,QZ)=\\[6pt]
2g((D^+_XQ)(Y),Z)-d^+F^+(X,J_+QY,J_+Z)-d^+F^+(X,J_+Y,J_+QZ)=\\[6pt]
d^+F^+(X,PY,Z)+d^+F^+(X,Y,PZ)\\[6pt]
\hspace{4cm}+2d^+F^+(X,J_-Y,J_+Z)+2d^+F^+(X,J_+Y,J_-Z).
\end{array}
\end{equation}
Therefore
$$
\begin{array}{l}
2\mathfrak{G} g((\nabla_{QX}Q)(Y),Z)=\mathfrak{G}[d^+F^+(QX,PY,Z)+d^+F^+(QX,Y,PZ)\\[6pt]
\hspace{3cm}+2d^+F^+(QX,J_-Y,J_+Z)+d^+F^+(QX,J_+Y,J_-Z)]
\end{array}
$$
Using the skew-symmetry of $d^{+}F^{+}$, it is easy to see that
$$
\begin{array}{l}
\mathfrak{G}[d^+F^+(QX,PY,Z)+d^+F^+(QX,Y,PZ)]=\\[6pt]
2\mathfrak{G}[d^{+}F^{+}(J_+J_-X,J_+J_-Y,Z)-d^{+}F^{+}(J_-J_+X,J_-J_+Y,Z)].
\end{array}
$$
We have $d^{+}F^{+}=-d^{-}F^{-}$, so  $d^{+}F^{+}$ is of type
$(2,1)+(1,2)$ for both $J_+$ and $J_-$. Therefore
$$
d^{+}F^{+}(A,B,C)=\mathfrak{G}d^{+}F^{+}(J_+A,J_+B,C)=\mathfrak{G}d^{+}F^{+}(J_-A,J_-B,C).
$$
It follows that
$$
\begin{array}{l}
\mathfrak{G} g((\nabla_{QX}Q)(Y),Z)=\\[6pt]
\hspace{2cm} \mathfrak{G}[d^{+}F^{+}(J_+J_-X,J_+J_-Y,Z)-d^{+}F^{+}(J_-J_+X,J_-J_+Y,Z)\\[6pt]
\hspace{2.4cm}+d^{+}F^{+}(J_+J_-X,J_-Y,J_+Z)-d^{+}F^{+}(J_-J_+X,J_-Y,J_+Z)\\[6pt]
\hspace{2.4cm}+d^{+}F^{+}(J_+J_-X,J_+Y,J_-Z)-d^{+}F^{+}(J_-J_+X,J_+Y,J_-Z)]=\\[8pt]
\hspace{2cm}  \mathfrak{G}[d^{+}F^{+}(J_-X,J_-Y,Z)-d^{+}F^{+}(J_+X,J_+Y,Z)\\[6pt]
\hspace{2.4cm}+d^{+}F^{+}(J_-X,J_-Y,Z)-d^{+}F^{+}(J_+X,Y,J_+Z)\\[6pt]
\hspace{2.4cm}+d^{+}F^{+}(J_-X,Y,J_-Z)-d^{+}F^{+}(J_+X,J_+Y,Z)]=\\[8pt]
\hspace{6cm}3[d^{+}F^{+}(X,Y,Z)-d^{+}F^{+}(X,Y,Z)]=0.
\end{array}
$$

This proves that $[Re\Pi,Re\Pi]=0$ which implies, as we have
mentioned, that $[\Pi,\Pi]=0$, i.e. $\Pi$ is a Poisson field.

One can also prove that the field $\Pi$ is Poisson using the fact
that the $2$-vector corresponding to the endomorphism
$J_{+}+J_{-}$ is Poisson \cite{LZ} and its $(2,0)$-part  is a
constant multiple of $\Pi$. {\it q. e. d.}

\vspace{.2in}



A holomorphic Poisson structure on a complex surface is merely a holomorphic section of its
anti-canonical bundle. Using this fact N. Hitchin \cite{Hit1} proposed a simple way for
constructing generalized K\"ahler structures on Del Pezzo surfaces. A different approach by M.
Gualtieri \cite{G3} based on the notion of generalized complex branes  extends this construction to
higher-dimensional Fano manifolds. Here, we state a modification of his result which can be proved
in the same way as \cite[Theorem 7.1]{G3}

\begin{thm}
Let $L$  be a holomorphic line bundle on an $n$-dimensional
compact complex manifold $M$ with holomorphic Poisson structure
$\sigma$ such that $c_1(L)^n\neq 0$. Let $(g_0,J_0)$ be a
pseudo-K\"ahler structure with K\"ahler form $F_0\in c_1(L)$.
Consider $\sigma$ and $F_0$ as homomorphisms $\sigma:(T^{\Bbb
C}M)^{\ast}\to T^{\Bbb C}M$ and $F_0:T^{\Bbb C}M\to (T^{\Bbb
C}M)^{\ast}$, and suppose that the following conditions are
satisfied:

(i) $\sigma\circ F_0=\overline{\partial}X^{1,0}$ for some $(1,0)$ vector field $X^{1,0}$;

(ii) $[Re\,X^{1,0}, Im\,\sigma]=0$ for the Schouten-Nijenhuis
bracket.

Then the choice of a Hermitian structure on $L$ with curvature $F_0$ determines a family of
generalized pseudo-K\"ahler structures $(g_t, J_t, J_0)$ with $J_t=\phi_t^*(J_0)$ for a 1-parameter
group of diffeomorphisms $\phi_t$ such that $J_t=J_0$ for $t\neq 0$ only at the poins of $M$ where
$\sigma=0$.

\end{thm}

\noindent {\bf Remark 1} Using Theorem 4 or the construction in \cite{Hit3}, one expects to produce
examples of generalized pseudo-K\"ahler structures on ruled surfaces over a Riemann surface of
genus greater than one. For example, consider a ruled surface $M$ over a curve $C$ of genus $g>1$
obtained as a projectivization of a vector bundle $V$ of degree $deg(V)<2-3g$. Its anti-canonical
bundle has a nowhere-vanishing  holomorphic section $s$ and the choice of a Hermitian metric on it
will produce a curvature 2-form $F_0=dd^clog |s|^2$. Suppose that $F_0$ is non-degenerate at each
point. Then Theorem 4 and \cite{Hit3} produce generalized pseudo-K\"ahler structure with
non-trivial canonical bundle. Note that when $V={\cal O}\oplus {\cal L}$ is decomposable, the
admissible metrics on $M$ considered in \cite{ACGT}  define Hermitian metrics on the anti-canonical
bundle of $M$ which are candidates to provide such $F_0$. However one can check that none of these
metrics  has a non-degenerate Ricci tensor. In case $deg(V)>2-2g$, there are metrics with this
property but there is no holomorphic Poisson structure. So, it is an open question whether any
ruled surface admits a generalized pseudo-K\"ahler structure. Note that  R. Goto \cite{Go} has
recently constructed positive definite generalized K\"ahler structures on some of these surfaces
using more general deformations of K\"ahler-Poisson structures \cite{Goto} than that considered in
\cite{G3}. However his approach is based on elliptic methods and can not be adapted directly to the
pseudo-Riemannian case.

\medskip

\section{Generalized pseudo-K\"ahler structures on 3-dimensional flag manifold}

 Consider the complex flag manifold $Fl=\{(L,V)|~ 0\in L\subset V\subset
\mathbb{C}^3, dim\,L=1,dim\,V=2\}$. It can be embedded into
$\mathbb{C}\mathbb{P}^2\times\mathbb{C}\mathbb{P}^{2}$ as the quadric
$Fl=\{(x_0,x_1,x_2;y_0,y_1,y_2)\in
\mathbb{C}\mathbb{P}^2\times\mathbb{C}\mathbb{P}^{2}|~x_0y_0+x_1y_1+x_2y_2=0\} $. Let $\omega$  be
the K\"ahler form of the standard K\"ahler structure on $\mathbb{C}\mathbb{P}^{2}$  normalized so
that $\omega$  be integral. Denote by $p_1$ and $p_2$ the projections of
$\mathbb{C}\mathbb{P}^2\times\mathbb{C}\mathbb{P}^{2}$ onto the first and the second factor. Set
$\omega_1=p_1^{\ast}\omega$ and $\omega_2=p_2^{\ast}\omega$. The restrictions of these forms to
$Fl$ will be denote by the same symbols.
\begin{lm}\label{nondegenerate}
For any integers $a$ and $b$ with $ab<0$
and $a+b\neq 0$, the form $F=a\omega_1+b\omega_2$ is
non-degenerate on $Fl$.
\end{lm}
{\it Proof:} Suppose that for such $a$ and $b$ the $2$-form
$F=a\omega_1+b\omega_2$ is degenerate at some point of $Fl$. The
group $U(3)$, embedded diagonally in $U(3)\times U(3)$, acts
transitively and holomorphically on $Fl$ and $F$ is invariant under this action. It
follows that $F$ degenerates at every point of $Fl$. This implies
that the top degree $F^3$ vanishes since $deg F=2$. We have
$\omega_i^3=(p_i^{\ast}\omega^3)|Fl=0$ for $i=1,2$. Therefore
$F^3=3ab(a\,\omega_1^2\wedge\omega_2+b\,\omega_1\wedge\omega_2^2)$.
Let $\psi:Fl\to Fl$ be the holomorphic map induced by the map
$\psi([x],[y])=([y],[x])$ on
$\mathbb{C}\mathbb{P}^2\times\mathbb{C}\mathbb{P}^{2}$. It  is
clear that $\psi^{\ast}\omega_1=\omega_2$ and
$\psi^{\ast}\omega_2=\omega_1$. Therefore
$0=\psi^{\ast}(a\omega_1^2\wedge\omega_2+b\omega_1\wedge\omega_2^2)=a\omega_1\wedge\omega_2^2+b\omega_1^2\wedge\omega_2$.
Then $(a+b)(\omega_1^2\wedge\omega_2+\omega_1\wedge\omega_2^2)=0$
and we get the identity
$(a+b)(\omega_1+\omega_2)^3=3(a+b)(\omega_1^2\wedge\omega_2+\omega_1\wedge\omega_2^2)=0$.
But the latter identity does not hold since $\omega_1+\omega_2$ is
the K\"ahler form of $Fl$ induced by the product of the
Fubini-Studi forms on each factor of $\mathbb{C}\mathbb{P}^2$, a
contradiction. {\it q. e. d.}

Later in the paper we'll need the following:

\begin{lm}\label{dbar}
Let $U$ and $V$ be commuting holomorphic vector fields on a
complex manifold and $\varphi$ a smooth function on the manifold.
Then $$(U\wedge V) \circ dd^c\varphi =
i\overline{\partial}((U\varphi)V-(V\varphi)U).$$
\end{lm}
{\it Proof}.  We use the identity
$d^c=\frac{1}{2}(i\overline{\partial}-i\partial)$ and the fact
that for any $(0,1)$-vector field $Z$,  $[U,Z]^{(1,0)}=0$. Then we have
$$
\begin{array}{lll}
2(dd^c\varphi)(U,Z)&=&iU(\overline{\partial}\varphi(Z))+iZ(\partial
\varphi(U))-i\overline{\partial}\varphi([U,Z])=\\[6pt]
& &iUZ\varphi+iZU\varphi-[U,Z]\varphi=2iZU\varphi,
\end{array}
$$
so
$$
\imath_Udd^c\varphi=i\overline{\partial}(U\varphi).
$$
From here we get
$$
\begin {array}{lll}
(U\wedge V)\circ dd^c\varphi = \imath_Udd^c\varphi\otimes V-\imath_Vdd^c\varphi\otimes U &=&
i\overline{\partial}(U\varphi)\otimes V -i\overline{\partial}(V\varphi)\otimes U=\\[6pt]
& &i\overline{\partial}((U\varphi)V-(V\varphi)U).
\end{array}
$$
{\it q. e. d.}

Now we are ready to prove the following:
\begin{prop}\label{flag}
The flag manifold $Fl$ admits a generalized pseudo-K\"ahler structure.
\end{prop}
{\it Proof:} Take arbitrary integers $a$ and $b$ with $ab<0$, $a+b\neq 0$. Then
$F_0=a\omega_1+b\omega_2$ is non-degenerate by Lemma \ref{nondegenerate}, so it determines a
pseudo-K\"ahler metric on $Fl$.

\smallskip

Since the  form $F_0$ is integral, it determines a Hermitian
holomorphic line bundle $L$ on $Fl$  with curvature $F_0$.
We have  $c_1(L)^3\neq 0$ since $c_1(L)^3$ is represented by the invariant form $F_0^3$ on $Fl$ and
the $2$- form $F_0$ is non-degenerate

\medskip
Now we want to define a holomorphic Poisson structure on $Fl$ as
$\sigma=Z_1\wedge Z_2$ for two commuting holomorphic vector fields
$Z_1$ and $Z_2$. Let $Z_1$ and $Z_2$ be the fields on
$\mathbb{C}\mathbb{P}^2\times\mathbb{C}\mathbb{P}^{2}$ generated
by the complex 1-parameter groups $(x_0, x_1, x_2; y_0, y_1, y_2)
\rightarrow (e^tx_0, e^{-t}x_1, x_2; e^{-t}y_0, e^ty_1, y_2)$ and
$(x_0, x_1, x_2; y_0, y_1, y_2)$ $ \rightarrow (e^{t}x_0, x_1,
e^{-t}x_2; e^{-t}y_0, y_1, e^{t}y_2)$, respectively. Clearly $Z_1$
and $Z_2$ are commuting holomorphic vector fields tangent to $Fl$.
Then $Z_1\wedge Z_2$ is a holomorphic Poisson structure on $Fl$.
To show that $Fl$ admits a generalized pseudo-K\"ahler structure
it remains only to check conditions $(i)$ and $(ii)$ in Theorem 4.
Denote by $X$ the holomorphic vector field on ${\Bbb C}{\Bbb P}^2$
generated by the group $(x_0, x_1, x_2) \rightarrow (e^tx_0,
e^{-t}x_1, x_2)$. Then $Z_1=(X\circ p_1,-X\circ p_2)$. Similarly
$Z_2=(Y\circ p_1,-Y\circ p_2)$ where $Y$ is the vector field on
${\Bbb C}{\Bbb P}^2$ generated by the group $(x_0, x_1, x_2)
\rightarrow (e^{t}x_0, x_1, e^{-t}x_2)$
 The bi-vector filed $\tau=X\wedge Y$ is a holomorphic section of the anti-canonical bundle of
${\Bbb C}{\Bbb P}^2$. Set $f=\ln ||\tau||^{2}$ where the norm is taken with respect to metric
yielded by the normalized Fubini-Study metric $g$ of ${\Bbb C}{\Bbb P}^2$. We claim that, although
$f$ is defined only outside of the zero set of $\tau$, the functions $Xf$ and $Yf$ are globally
defined and smooth. To check this we use the standard coordinates of $\mathbb{C}P^2$. For the
coordinates $z_1=\displaystyle{\frac{x_1}{x_0}}$, $z_2=\displaystyle{\frac{x_2}{x_0}}$, set
$$
g_{\alpha\overline\beta}=g(\frac{\partial}{\partial
z_{\alpha}},\frac{\partial}{\partial \overline z_{\beta}})\quad
\mbox{and}\quad G_{(z)}=g_{1\overline 1}g_{2\overline
2}-|g_{1\overline 2}|^2.
$$
Then $||\tau||^2=4|z_1z_2|^2G_{(z)}$ and we have
\begin{equation}\label{z}
\begin{array}{l}
\displaystyle{X=-2z_1\frac{\partial}{\partial
z_1}-z_2\frac{\partial}{\partial z_2},\quad
Y=-z_2\frac{\partial}{\partial z_2},\quad
\tau=2z_1z_2\frac{\partial}{\partial
z_1}\wedge\frac{\partial}{\partial z_2}},\\[10pt]
\displaystyle{Xf=-3-2z_1\frac{\ln G_{(z)}}{\partial
z_1}-z_2\frac{\ln G_{(z)}}{\partial z_2},\quad Yf=-1-z_2\frac{\ln
G_{(z)}}{\partial z_2}}.
\end{array}
\end{equation}

In the coordinates $u_1=\displaystyle{\frac{x_0}{x_1}}$,
$u_2=\displaystyle{\frac{x_2}{x_1}}$ we have
\begin{equation}\label{u}
\begin{array}{l}
\displaystyle{X=2u_1\frac{\partial}{\partial
u_1}+u_2\frac{\partial}{\partial u_2},\quad
Y=-u_2\frac{\partial}{\partial u_2},\quad
\tau=-2u_1u_2\frac{\partial}{\partial
z_1}\wedge\frac{\partial}{\partial z_2}},\\[10pt]
\displaystyle{ Xf=3+2u_1\frac{\ln G_{(u)}}{\partial
u_1}+u_2\frac{\ln G_{(u)}}{\partial u_2},\quad Yf=-1-u_2\frac{\ln
G_{(u)}}{\partial u_2}}.
\end{array}
\end{equation}
Finally, in the coordinates $v_1=\displaystyle{\frac{x_0}{x_2}}$,
$v_2=\displaystyle{\frac{x_1}{x_2}}$ we have
\begin{equation}\label{v}
\begin{array}{l}
\displaystyle{X=v_1\frac{\partial}{\partial
v_1}-v_2\frac{\partial}{\partial v_2},\quad
Y=v_1\frac{\partial}{\partial v_1}+v_2\frac{\partial}{\partial
v_2},\quad
\tau=2v_1v_2\frac{\partial}{\partial v_1}\wedge\frac{\partial}{\partial v_2}},\\[10pt]
\displaystyle{Xf=v_1\frac{\ln G_{(v)}}{\partial v_1}-v_2\frac{\ln
G_{(v)}}{\partial v_2},\quad Yf=2+v_1\frac{\ln G_{(v)}}{\partial
v_1}-v_2\frac{\ln G_{(v)}}{\partial v_2}}.
\end{array}
\end{equation}
It follows from (\ref{z}), (\ref{u}), (\ref{v}) that $\tau$
vanishes on the analytic set $C=\{[x]\in {\Bbb C}{\Bbb P}^2:
x_0x_1x_2=0\}$ and that $Xf$, $Yf$ can be extended to smooth
functions on a neighborhood of every point of $C$. Since ${\Bbb
C}{\Bbb P}^2\setminus C$ is dense, we see that $Xf$, $Yf$ can be
extended to  unique smooth functions on the whole space ${\Bbb
C}{\Bbb P}^2$. We shall denote the extensions by the same symbols.
Identities (\ref{z}), (\ref{u}), (\ref{v}) imply also that if
$\zeta=(\zeta_1,\zeta_2)$ is a standard coordinate system of
${\Bbb C}{\Bbb P}^2$, we have $dd^{c}\ln||\tau||^2=dd^{c}\ln
G_{(\zeta)}$ on ${\Bbb C}{\Bbb P}^2\setminus C$. Therefore
$dd^{c}\ln||\tau||^2$ on ${\Bbb C}{\Bbb P}^2\setminus C$ is the
Ricci from of the standard K\"ahler structure on ${\Bbb C}{\Bbb
P}^2$. As it is well-known, the Ricci form of this structure is
equal to $3$ times the K\"ahler form. Thus, since we are working
with the normalized K\"ahler form, we have
$dd^{c}\ln||\tau||^2=3\lambda\omega$ where $\lambda>0$ is a
constant. Hence, for $k=1,2$, $dd^{c}(\ln||\tau||^2\circ
p_k)=3\lambda p_k^{\ast}\omega=3\lambda\,\omega_k$ on the set
$M=\{(x_0,x_1,x_2;y_0,y_1,y_2)\in
\mathbb{C}\mathbb{P}^2\times\mathbb{C}\mathbb{P}^{2}|~x_0x_1x_2y_0y_1y_2\neq
0\} $. Thus on $M$ we have
\begin{equation}\label{circ}
(Z_1\wedge Z_2)\circ F_0=\frac{1}{3\lambda}(Z_1\wedge Z_2)(a\,
dd^{c}(\ln||\tau||^2\circ p_1)+b\, dd^{c}(\ln||\tau||^2\circ p_2))
\end{equation}

It follows from (\ref{circ}) and Lemma~\ref{dbar} that if we set
$$
X^{1,0}=\frac{i}{3\lambda}\{[a(Xf)\circ p_1-b(Xf)\circ
p_2]Z_2-[a(Yf)\circ p_1-b(Yf)\circ p_2]Z_1\}
$$
where $f=\ln ||\tau||^{2}$ as above, we have $(Z_1\wedge Z_2)\circ F_0=\overline{\partial}X^{1,0}$
on the open set $M$. This identity holds everywhere since the vector field $X^{1,0}$ is smooth on
$\mathbb{C}\mathbb{P}^2\times\mathbb{C}\mathbb{P}^{2}$  and $M$ is dense.  Thus condition $(i)$ of
Theorem 4 is satisfied for $\sigma=Z_1\wedge Z_2$. To show that condition $(ii)$ also holds, we
note that
$$
[X^{1,0},Z_1\wedge Z_2]=-\frac{i}{3\lambda}\{a([X,Y]f)\circ
p_1+b([X,Y]f)\circ p_2\}Z_1\wedge Z_2=0
$$
since $[X,Y]=0$. The function $f$ is real-valued, so
$\overline{Xf}=\overline{X}f$, $\overline{Yf}=\overline{Y}f$ and
we have
$$
[\overline{X^{1,0}},\overline Z_1\wedge \overline
Z_2]=\frac{i}{3\lambda}\{a(\overline{[X,Y]}f)\circ
p_1+b(\overline{[X,Y]}f)\circ p_2\}\overline Z_1\wedge \overline
Z_2=0.
$$
Using the identities $[X,\overline X]=[Y,\overline Y]=[X,\overline
Y]=[\overline X,Y]=0$, it is easy to see that
$$
[X^{1,0},\overline Z_1\wedge \overline
Z_2]-[\overline{X^{1,0}},Z_1\wedge Z_2]=0.
$$
It follows that $[Re\,X^{1,0},Im(Z_1\wedge Z_2)]=0$. Then, by Theorem 4, the flag manifold $Fl$
admits a generalized pseudo-K\"ahler structure. {\it q. e. d.}

Note that $Fl$ admits also  a usual
generalized K\"ahler structure \cite{Goto}.

\begin{cor}
Any holomorphic line bundle on the 3-dimensional flag manifold $Fl$ carries a structure of a
holomorphic Poisson module with respect to the holomorphic Poisson structure $U_1\wedge U_2$
defined by commuting holomorphic vector fields $U_1$ and $U_2$.

\end{cor}
{\it Proof:} First we notice that any two commuting vector fields on $Fl$ span a maximal torus in
the algebra $sl(3,\mathbb{C})$ of the holomorphic vector fields on $Fl$ and all such tori are
conjugate in the group of biholomorphisms. So, we may assume that the vector fields $U_1$ and $U_2$
in the corollary coincide $Z_1$ and $Z_2$ defined in the proof of Proposition \ref{flag}. Denote by
$K$ the canonical bundle of $\mathbb{C}\mathbb{P}^2$. It is well known that every holomorphic line
bundle over $Fl$ is of the form
$L_{mn}=\displaystyle{\frac{m}{3}p_1^{\ast}K+\frac{n}{3}p_2^{\ast}K}$ where $m,n\in \mathbb{Z}$. If
we consider $K$ with the metric induced by the normalized Fubini-Study metric of
$\mathbb{C}\mathbb{P}^2$, the curvature form of $K$ with respect to its canonical connection is
equal to the K\"ahler form $\omega$. Therefore the form
$\displaystyle{\frac{m}{3}p_1^{\ast}\omega+\frac{n}{3}p_2^{\ast}\omega=\frac{m}{3}\omega_1+\frac{n}{3}\omega_2}$
represents the first Chern class of $L_{mn}$. Denote this form by $F$ and set $\sigma=Z_1\wedge
Z_2$. We have seen above that there is a $(1,0)$-vector field $X^{1,0}$ such that $\sigma\circ
X^{1,0}=\overline{\partial} X^{1,0}$ and $[X^{1,0},\sigma]=0$. Now the Corollary follows from
\cite[Proposition 10]{G3} since the first Chern class of $L_{mn}$ coincides with its Atiyah class.
{\it q. e. d.}

\section{The four-dimensional case}

In dimension four, a pseudo-hermitian metric is either positive
(negative) definite or of signature (2,2). Using the results in
Section 4 we shall prove the following:

\begin{thm}
Let $(M,g,J_+,J_-)$ be a compact  pseudo-bihermitian 4-manifold.

\smallskip

$(i)$ If $d^+F^+ = -d^-F^- $, then $(M, J_+)$ (and $(M,J_-)$) is
one of the following complex surfaces:  a complex torus, a K3
surface, a primary Kodaira surface, a blow-up of a surface of
class $VII_0$, a ruled surface described in \cite{B} with $\chi
{\pm} \tau$ divisible by 4, where $\chi$ and $\tau$ are the Euler
characteristic and the signature of $M$.

\smallskip

$(ii)$ If the bihermitian structure is strict, then $(M, J_+)$
(and $(M,J_-)$) is one of the following: a complex torus, a K3
surface, a primary Kodaira surface, a properly elliptic surface of
odd first Betti number, a Hopf surface, a minimal Inoue surface
without curves.
\end{thm}
{\it Proof:} According to Theorem~\ref{Poisson}, under assumption $(i)$ there is a non-zero
holomorphic section of the anti-canonical bundle of $(M,J_+)$. Such surfaces with even first Betti
number are described in \cite{B} and they exhaust the first four cases in $(i)$. The restriction on
$\chi {\pm} \tau$ in the last case comes from Matsushita's topological condition for existence of a
split-signature metric \cite{Mat}. For the case of surfaces with odd first Betti number, we notice
that the proof of Proposition 2.3 in \cite{Bo} shows that either the Kodaira dimension of $(M,J_+)$
(and $(M,J_-)$) is $-\infty$ or its canonical bundle is holomorphically trivial. Then the Kodaira
classification of minimal compact complex surfaces \cite{BPV} leads to the list in $(i)$.

Part $(ii)$ follows from the fact that the canonical bundle is
topologically trivial in the case of strictly pseudo-bihermitian
surfaces, since the $2$-form $\Omega^{(0,2)}$ given by (\ref
{omega}) provides a non-vanishing section. So one can use the
well-known list of the surfaces with vanishing first Chern class
\cite{W} 
 {\it q. e. d.}

\vspace{.2in}

\noindent {\bf Remark} {\bf 2}.Notice that, by \cite[Lemma 2.1]{B}, if a compact complex surface is
not minimal and has a nowhere-vanishing holomorphic section of the anti-canonical bundle, then its
minimal model also admits such a section. Moreover, the dimension of the space of holomorphic
sections decreases by at most one after a blow-up. It keeps the same dimension only if the blow-up
is at a base point of the anti-canonical linear system. This leads to additional restrictions on
the possible blow-ups of surfaces in case $(i)$, but we shall not discuss this question here.

\smallskip

\noindent {\bf Remark 3}. There are generalized pseudo-K\"ahler manifolds $(M,g, J_+, J_-)$ so that
$J_+$ and $J_-$ induce opposite orientations. In the four dimensional case such structures commute.
In any dimension, for a generalized pseudo-K\"ahler manifold with commuting  $J_+$ and $J_-$, the
same reasoning as in \cite{AGu} shows that the holomorphic tangent bundle of $(M,J_+)$ splits into
a sum of two holomorphic subbundles. Conversely, if the holomorphic tangent bundle of a compact
complex surface $(M,J)$ splits, then by \cite{AGu} there is a generalized (pseudo) K\"ahler
structure $(g,J_+,J_-)$ such that $J_+=J$ and $[J_+,J_-]=0$.

\medskip

\section{ Generalized pseudo-K\"ahler structures via deformations of para-hyperk\"ahler structures}
\medskip

 It has been observed in \cite{AGG,G1,Hit1} that one can explicitly
define a generalized K\"ahler structure by means of a hyperk\"ahler structure. Given a
para-hyperk\"ahler structure, a similar construction can be applied to obtain a generalized
pseudo-K\"ahler structure. Let $\{g,J_1, J_2, J_3\}$ be a para-hyperk\"ahler structure on a
$4$-manifold $M$ with $J_1^2=-J_2^2=-J_3^2=-Id$ and $J_3=J_1J_2$. We would like to construct two
commuting generalized almost complex structures $I_1$ and $I_2$ following \cite{Hit1}. To do this
we need complex valued $2$-forms $\beta_1$ and $\beta_2$ on $M$ which satisfy
\begin{equation}\label{1}
(\beta_1-\beta_2)^2=(\beta_1-\overline{\beta_2})^2=0,~
\beta_1\neq\beta_2,~ \beta_1\neq\overline{\beta_2}
\end{equation}
at every point. We set
$\exp(\beta_k)=1+\beta_k+\frac{1}{2}\beta_k^2$, $k=1,2$, and
$(X+\xi).\exp(\beta_k) = \imath_X
\exp(\beta_k)+\xi\wedge\exp(\beta_k)$ for $X+\xi\in TM\oplus
T^{\ast}M$ (the Clifford action of $TM\oplus T^{\ast}M$  on the
forms). Then $E_k=\{A\in (TM\oplus T^{\ast}M)^{\Bbb
C}~|~A.\exp(\beta_k)=0\}$ is the $+i$-eigenspace of a generalized
almost complex structure $I_k$. If $\beta_k$ is closed, $I_k$ is
Courant integrable \cite{G1, Hit1}. It is shown in \cite[Lemma
1]{Hit1} that $I_1$ and $I_2$ commute. Moreover, $E_1\cap
E_2\oplus\overline{E_1}\cap\overline{E_2}$ is the
$(-1)$-eigenspace of $I_1I_2$ and $E_1\cap\overline{E_2}\oplus
\overline{E_1}\cap E_2$ is the $(+1)$-eigenspace. Note also that
$E_1\cap E_2=\{U-\imath_U\beta_1~|~U\in T^{\Bbb C}M,~
\imath_U\beta_1=\imath_U\beta_2\}$ (\cite{Hit1}). Thus, for
$A=U-\imath_U\beta_1\in E_1\cap E_2$, $B=V-\imath_V\beta_1\in
E_1\cap E_2$, we have
\begin{equation}\label{metric}
\begin{array}{c}
<A+\overline{A},B+\overline{B}>=-Re\{(\beta_1-\overline{\beta_1})(U,\overline{V})\}=-Re\{(\beta_2-\overline{\beta_2})(U,\overline{V})\}\\[6pt]
=-Re\{(\beta_1-\overline{\beta_2})(U,\overline{V})\}
\end{array}
\end{equation}
Now, given a para-hyperk\"ahler structure  $\{g,J_1, J_2, J_3\}$ on a $4$-manifold $M$, set
$J_{+}=J_1$ and $J_{-}=aJ_1+bJ_2+cJ_3$ where $a,b,c$ are fixed numbers such that $a^2-b^2-c^2=1$
and $a\neq 1$. Then $J_{+}$ and $J_{-}$ are complex structures compatible with the metric $g$
satisfying the identity
\begin{equation}\label{-a}
J_{+}J_{-}+J_{-}J_{+}=-2aId.
\end{equation}
As in Section 2, set
$$
K=\displaystyle{\frac{1}{2\sqrt{a^2-1}}}[J_+,J_-],~ S_{+}=-\displaystyle{\frac{1}{\sqrt{a^2-1}}
}(J_- -aJ_+)
$$
Then $\{g,J_{+},K,S_{+}\}$ is a para-hyperhermitian structure with $S_{+}=J_{+}K$. Let
$F^{+}(X,Y)=g(J_{+}X,Y)$, $F^{K}(X,Y)=g(KX,Y)$ and $\omega'(X,Y)=g(S_{+}X,Y)$ be the corresponding
fundamental $2$-forms. Similarly, if
$$
S_{-}=\displaystyle{\frac{1}{\sqrt{a^2-1}} }(J_{+} -aJ_{-}),
$$
then $\{g,J_{-},K,S_{-}\}$ is a para-hyperhermitian structure with $S_{-}=J_{-}K$. We denote the
fundamental $2$-forms of $J_{-}$ and $S_{-}$ by $F^{-}$ and $\omega''$, respectively. Set
$$
\omega_{+}=\omega' + \omega'',\quad \omega_{-}=\omega' - \omega''.
$$
Then
\begin{equation}\label{pm}
\begin{array}{l}
\omega_{+}(X,Y)=\displaystyle{\sqrt{\frac{a+1}{a-1}}g(J_{+}X-J_{-}X,Y)=\sqrt{\frac{a+1}{a-1}}(F^{+}(X,Y)-F^{-}(X,Y))}\\[8pt]
\omega_{-}(X,Y)=\displaystyle{\sqrt{\frac{a-1}{a+1}}g(J_{+}X+J_{-}X,Y)=\sqrt{\frac{a-1}{a+1}}(F^{+}(X,Y)+F^{-}(X,Y))}\\[8pt]
\end{array}
\end{equation}
In particular, the forms $\omega_{+}$ and $\omega_{-}$ are closed since $F^{+}$ and $F^{-}$ are so.
Identity (\ref{nablaQ}) implies that $\nabla [J_{+},J_{-}]=0$, thus $\nabla K=0$. Therefore the
form $F^{K}$ is also closed. Now, similar to \cite{G1} we set
$$
\beta_1=F^K+i\omega_{+}, \quad \beta_2=-F^K+i\omega_{-}.
$$
Conditions (\ref{1}) for these forms are equivalent to
\begin{equation}\label{bih}
F^K\omega_+=F^K\omega_-=\omega_+\omega_-=\omega_+^2+\omega_-^2-4(F^K)^2=0.
\end{equation}
Let $X$ be a tangent vector with $g(X,X)=1$. Then
$\{X,J_{+}X,KX,S_{+}X\}$ is a $g$-orthonormal basis of tangent
vectors. Using (\ref{-a}), (\ref{pm}) and the paraquaternionic
identities, it is easy to see that

$$
\begin{array}{ccl}
(\omega_{+}\wedge\omega_{+})(X,J_{+}X,KX,S_{+}X)=4(a+1),\\[6pt]
(\omega_{-}\wedge\omega_{-})(X,J_{+}X,KX,S_{+}X)=-4(a-1),\\[6pt]
(\omega_{+}\wedge\omega_{-})(X,J_{+}X,KX,S_{+}X)=0,\quad
(F^K\wedge\omega_{\pm})(X,J_{+}X,KX,S_{+}X)=0.
\end{array}
$$
We also have $(F^K\wedge F^K)(X,J_{+}X,KX,S_{+}X)=2$. It follows that identities (\ref{bih}) are
satisfied.

The identity $\beta_1-\beta_2=2F^K+i(\omega_{+}-\omega_{-})$
implies that a vector $U\in T^{\Bbb C}M$ satisfies
$\imath_U(\beta_1-\beta_2)=0$ if and only if
\begin{equation}\label{imath}
\sqrt{a^2-1}KU+iJ_{+}U-iaJ_{-}U=0.
\end{equation}
Thus $E_1\cap E_2=\{U-\imath_U\beta_1~|~U\in T^{\Bbb C}M,~ U~ \mbox {satisfies} ~(\ref{imath})\}$.
Let $L^{-}$ be the $(-1)$-eigenspace of $I_1I_2$ acting on $TM\oplus T^{\ast}M$. Any $X+\xi\in
L^{-}$ can be written as $X+\xi=U+\overline U$ where $U=\frac{1}{2}(X+iY)\in E_1\cap E_2$, $Y\in
TM$, and $\xi=\imath_U\beta_1-\imath_{\overline U}\overline{\beta_1}$. In this notation,
(\ref{imath}) is equivalent to
\begin{equation}\label{ima}
\begin{array}{ll}
\sqrt{a^2-1}KX-J_{+}Y+aJ_{-}Y=0\\[6pt]
\sqrt{a^2-1}KY+J_{+}X-aJ_{-}X=0.
\end{array}
\end{equation}
In fact, either of these identities is a consequence of the other one. For every
$V=\frac{1}{2}(Z+iT)\in E_1\cap E_2$, we have
$$
\begin{array}{c}
Re{(\beta_1-\overline\beta_1)(U,\overline V)}=\displaystyle{\sqrt{\frac{a+1}{a-1}}[g(J_{+}X-J_{-}X,T)-g(J_{+}Y-J_{-}Y,Z)]}= \\[8pt]
-\displaystyle{\sqrt{\frac{a+1}{a-1}}[g(X,J_{+}T-J_{-}T)+g(J_{+}Y-J_{-}Y,Z)]}.
\end{array}
$$
Applying $K$ to the second identity of (\ref{ima}) we get $\sqrt{a^2-1}Y=S_{+}X-aS_{-}X$. This
gives
$$
\begin{array}{l}
\sqrt{a^2-1}J_{+}Y=-KX+\displaystyle{\frac{a}{\sqrt{a^2-1}}}X+\displaystyle{\frac{a^2}{\sqrt{a^2-1}}}J_{+}J_{-}X,\\[8pt]
\sqrt{a^2-1}J_{-}Y=aKX+\displaystyle{\frac{1}{\sqrt{a^2-1}}}X+\displaystyle{\frac{a}{\sqrt{a^2-1}}}J_{-}J_{+}X.
\end{array}
$$
It follows that
$$
\sqrt{a^2-1}(J_{+}Y-J_{-}Y)=(a-1)(a+2)KX+\displaystyle{\frac{a-1}{\sqrt{a^2-1}}}X.
$$
Similarly,
$$
\sqrt{a^2-1}(J_{+}T-J_{-}T)=(a-1)(a+2)KZ+\displaystyle{\frac{a-1}{\sqrt{a^2-1}}}Z.
$$
Then
\begin{equation}\label{ff}
(a-1)Re{(\beta_1-\overline\beta_1)(U,\overline V)}=-\displaystyle{\sqrt{\frac{a-1}{a+1}}}g(X,Z).
\end{equation}

Suppose that $<X+\xi,A>=0$ for every $A\in L^{-}$. Take any $Z\in TM$ and set
$$
T=(a^2-1)^{-1/2}[S_{+}Z-aS_{-}Z].
$$
Then $V=\frac{1}{2}(Z+iT)$ satisfies (\ref{imath}). Indeed we have
$$
\begin{array}{l}
\sqrt{a^2-1}KZ-J_{+}T+aJ_{-}T=\\[6pt]
\displaystyle{\sqrt{a^2-1}KZ-\frac{1}{\sqrt{a^2-1}}(-KZ-aJ_{+}S_{-}Z)+\frac{a}{\sqrt{a^2-1}}(J_{-}S_{+}Z+aKZ)}= \\[6pt]
\displaystyle{\frac{1}{\sqrt{a^2-1}}(2a^2KZ+a(J_{+}S_{-}Z+J_{-}S_{+}Z))=\frac{1}{\sqrt{a^2-1}}(2a^2KZ-a^2\frac{[J_{+},J_{-}]Z}{\sqrt{a^2-1}})}=\\[6pt]
\displaystyle{\frac{1}{\sqrt{a^2-1}}}(2a^2KZ-2a^2KZ)=0.
\end{array}
$$
Moreover,
$$
\sqrt{a^2-1}KT+J_{+}Z-aJ_{-}Z=KS_{+}Z-aKS_{-}Z+J_{+}Z-aJ_{-}Z=0.
$$
Thus $V\in E_1\cap E_2$ and, by our assumption, (\ref{metric}) and (\ref{ff}), we have $g(X,Z)=0$.
Since the latter identity holds for every $Z$, we conclude that $X=0$. Then
$Y=(a^2-1)^{-1/2}[S_{+}X-aS_{-}X]=0$, hence $U=0$, thus $\xi=\imath_U\beta_1-\imath_{\overline
U}\overline{\beta_1}=0$. This proves that the canonical inner product on $TM\oplus T^*M$ is
non-degenerate on $L^{-}$. Moreover, the inclusion $TM\cap L^{-}\subset E_1\cap E_2$ and identity
(\ref{imath}) imply that $TM\cap L^{-}=\{0\}$. Similar arguments show that the metric $<.~,~.>$ is
non-degenerate on the $(+1)$-eigenspace $L^{+}$ of $I_1I_2$ and $TM\cap L^{+}=\{0\}$. Thus
$I_1,I_2$ is a generalized pseudo-K\"ahler structure on $M$.

We can deform this structure using arbitrary smooth function $f$ on $M$. Let $H_t$ be the flow of
the $F^K$-Hamiltonian vector field $\imath_{df}F^K$, so $H_t^*(F^K)=F^K$. Define
$$
\gamma_1=F^K+i(\omega'+H_t^*\omega''), \quad \gamma_2=-F^K+i(\omega'-H_t^*\omega'').
$$
Then $\gamma_1-\gamma_2=2F^K+2iH_t^*\omega''=H_t^*(2F^K+2i\omega'')=H_t^*(\beta_1-\beta_2)$ and
$\gamma_1-\overline\gamma_2=\beta_1-\overline\beta_2$. It follows that for small $t$, the forms
$\gamma_1$ and $\gamma_2$ define a generalized pseudo-K\"ahler structure .

  Finally, let us note that a generalized pseudo-K\"ahler structure can be explicitly defined by
means of the pseudo-K\"ahler structures $(g,J_{+})$, $(g,J_{-})$ and \cite[(6.14)]{G1}.

\noindent {\bf Example 1.} The construction above can be applied to 4-tori and primary Kodaira
surfaces since each of these surfaces admits a para-hyperk\"ahler structure (see, for example,
\cite{Ka,Kam}). Recall that the Kodaira surfaces do not admit any (positive) generalized K\"ahler
structure \cite{A, AGu}.

\medskip

\noindent {\bf Example 2.}  Any para-hyperhermitian structure which is locally conformally
para-hyperk\"ahler can be deformed as in \cite{AGG} to obtain a strictly pseudo-bihermitian
structure. The universal cover of the locally conformally para-hyperk\"ahler manifold $M$ is
globally conformally para-hyperk\"ahler. The deformation is performed on its para-hyperk\"ahler
structure such that $H_t$ is invariant with respect to the fundamental group of $M$. Then one
obtains a generalized pseudo-K\"ahler structure which after a (global) conformal change descends to
a pseudo-bihermitian structure on the quotient. In particular, there are pseudo-bihermitian metrics
on properly elliptic surfaces of odd first Betti number and the Inoue surfaces of type $S^+$
\cite{DGMY2}. These surfaces do not admit any (positive) bihermitian structure \cite{A}.  On the other hand the 
quaternionc Hopf surfaces admit both bihermitian and
pseudo-bihermitian structures since they have both hyperhermitian and para-hyperhermitian metrics
\cite{DGMY2}. They also have bihermitian metrics arising from twisted generalized K\"ahler
structures \cite{AGu}, however it is not
clear whether these surfaces admit twisted generalized pseudo-K\"ahler structures. The same
question is open for K3 surfaces too.

Notice that the above constructions produce "complementary" examples of bihermitian and
pseudo-bihermitian structures  on
the surfaces in the lists in Theorem 9. We summarize the examples
obtained so far in:

\begin{prop} Generalized pseudo-K\"ahler structures exist on complex 2-tori and primary Kodaira surfaces. Pseudo-bihermitian 
structures exist also on the quaternionic Hopf
surfaces, properly elliptic surfaces with odd first Betti number
and Inoue surfaces of type $S^+$.
\end{prop}

\noindent{\bf Example 3.} Here we provide examples of complex structures $J_+$ and $J_-$ satisfying
the relation (\ref{anticom})
\[
J_+J_-+J_-J_+=2pId
\]
for a nonconstant function $p$ with $|p|>1$, which are not compatible with any global neutral
metric. Consider Example 1 above in the case of a complex torus which is a product of 2 elliptic
curves. It admits a holomorphic involution $\phi$ without fixed points, such that the quotient is a
smooth complex surface. This surface is called a hyperelliptic surface of type $I_a$. One can check
that the natural para-hypercomplex structure of the torus descends to a para-hypercomplex structure
on the quotient, but it admits no compatible para-hyperhermitian metrics \cite{DGMY2}. In
particular, one can fix a para-hyperk\"ahler family of $\phi$-invariant complex structures on the
torus and can deform any two structures of this family via the procedure described in Example 2.
The Hamiltonian deformations $H_t$ are defined by a single function and if one chooses this
function to be $\phi$-invariant, then both $(J_+)_t=J_+$ and $(J_-)_t$ are $\phi$-invariant for all
$t$. Since they satisfy the relation (\ref{anticom}) for small $t$, they descend to structures
which satisfy the same identity on the quotient hyperelliptic surface. Since $|p|>1$ at any point
for fixed $t$, $K \neq 0$ everywhere. If there were a compatible metric, then the fundamental forms
$F^K+iF^{J_+K}$ obtained as in the consideration above would provide a trivialization of the
canonical bundle, which is an absurd because the canonical bundle of a hyperelliptic surface is not
topologically trivial.
\medskip

\section{Null-planes of 4-dimensional pseudo-biher-mitian
metrics}

In this section we show that, under a mild restriction, a naturally defined null-plane distribution
on a pseudo-bihermitian 4-manifold $M$ determines a local Engel structure. Recall that an {\bf
Engel structure} is by definition a 2-dimensional distribution ${\cal D}$ on a 4-manifold $M$ such
that $rank[{\cal D},{\cal D}]=3$ and $rank[{\cal D},[{\cal D},{\cal D}]]=4$ at each point of $M$.
These structures have been actively investigated recently (see the introduction in \cite{V} for an
overview). They admit
 canonical coordinates and are preserved by small ${\cal
 C}^2$-deformations. The global existence of an oriented Engel structure on an oriented
compact manifold leads to triviality of its tangent bundle.
Moreover, Vogel \cite{V} showed  that the converse also holds -
any paralellizable  4-manifold admits such a structure.

Let $(M,g, J_+, J_-)$ be a pseudo-bihermitian 4-manifold with $J_+J_-+J_-J_+=2pId$ where $|p|>1$.
Let $F^{\pm}$ and $\theta_{\pm}$ be the K\"aher and the Lee form of $(g,J_{\pm})$, respectively.
Suppose that the pseudo-bihermitian structure is defined by a (twisted) generalized pseudo-K\"ahler
one. Then $d^+F^+ + d^-F^-=0$ by Theorem~\ref{gK} and taking the Hodge-dual 1-forms we get
$\theta_+ +\theta_-=0$.

If we set $K =[J_+,J_-]/2\sqrt{p^2-1}$ as above, then $K^2=Id$ and
$K\neq \pm Id$. Moreover, $g(KX,Y)=-g(X,KY)$, in particular the
eigenspaces of $K$ consists of isotropic vectors.
\begin{lm}\label{N}
For the  endomorphism $N_{\pm}=J_{+}+(p\pm\sqrt{p^2-1})J_{-}$ of
$TM$, we have $Ker\,N_{\pm}=Im\,\,N_{\pm}=\mbox{$\mp$ 1-eigenspace
of $K$}$.
\end{lm}
{\it Proof:}  It is easy to see that $N_{\pm}^2=0$ and $Ker\,N_{+}\cap Ker\,N_{-}=\{0\}$. Moreover,
\begin{equation}\label{JK}
\begin{array}{l}
-KJ_{+}=J_{+}K=\displaystyle{\frac{p}{\sqrt{p^2-1}}J_{+}-\frac{1}{\sqrt{p^2-1}}J_{-}},\\[8pt]
-KJ_{-}=J_{-}K=\displaystyle{\frac{1}{\sqrt{p^2-1}}J_{+}-\frac{p}{\sqrt{p^2-1}}J_{-}}.
\end{array}
\end{equation}
It follows that if ${\cal K}_{\pm}$ is the $\pm 1$-eigenspace of
$K$, then ${\cal K}_{-}\subset Ker\, N_{+}$ and ${\cal
K}_{+}\subset Ker\, N_{-}$. Hence $dim\,Ker\, N_{\pm}\geq 2$. This
implies the lemma since the kernels of $N_{+}$ and $N_{-}$  are
transversal and the dimension of the ambient space is $4$. {\it q.
e. d.}
\smallskip

Denote the vector field dual to $\theta_{\pm}$ w.r.t. $g$ by the same letter. Set
$X=(J_{+}+fJ_-)\theta_+$ where  $f=p-\sqrt{p^2-1}$ and $Y=\theta_{+}+K\theta_+$. Clearly $X,Y\in
Ker N_-$. One can easily see that $X$ and $Y$ are isotropic. Assume that $|\theta_+|_x\neq 0$ at
some point $x\in M$. Then the vector fields $X$ and $Y$ are linearly independent at $x$. Indeed,
suppose that $\lambda X +\mu Y=0$ at $x$ for some constants $\lambda$ and $\mu$. Applying $N_{-}$
to both sides of this identity, we get $\mu(N_- +N_-K)\theta_{+}=0$ at $x$. Then, using (\ref{JK}),
we compute easily that $\mu(J_{+}\theta_{+}-fJ_{-}\theta_{+})_x=0$. If
$J_{+}\theta_{+}=fJ_{-}\theta_{+}$ at $x$, we would have $|\theta_{+}|_x=f(x)|\theta_{+}|_x$, hence
$|\theta_{+}|_x=0$, a contradiction. Therefore $\mu=0$, thus $\lambda(J_+\theta_+
+fJ_-\theta_+)_x=0$. This implies $\lambda=0$ since $|\theta_+|_x\neq 0$.

Now define a $2$-plane in $T_xM$ setting
\begin{equation}\label{D}
{\cal D}_x=Span(X,Y)_x.
\end{equation}

\begin{thm}Let $(M,g,J_+,J_-)$ be a (twisted) generalized pseudo-K\"ahler
$4$-manifold with nowhere vanishing Lee forms $\theta_{+} = -\theta_-$ and such that $J_{+}J_{-}
+J_{-}J_{+}=2p Id$  with $|p|>1$.  Then  the null distribution ${\cal D}$ defined by (\ref{D}) is
an Engel structure on an open subset of $M$ or the flow of $Y$ consists of null-geodesics.
\end{thm}
{\it Proof:} Set $N=N_{-}$. Then ${\cal D}=\,Ker\, N = \,Im\, N$. We are going to calculate
$N[X,Y]$ and show that it is proportional to $NJ_+Y$. This will imply that $[X,Y]\in
Span(X,Y,J_+Y)$, so $rank[{\cal D},{\cal D}]=3$. Then we will show that $[Y,J_+Y]$ has vanishing
$J_+X$ component iff $\nabla_YY= FY$ for some smooth function $F$. This proves that either the flow
of $Y$ is geodesic or $rank[{\cal D},[{\cal D},{\cal D}]]=4$ on an open subset of $M$.

For the Levi-Civita connection we
have \cite{AGG}:
$$2(\nabla_X
J_{\pm})Y=g(X,Y)J_{\pm}\theta_{\pm}+g(J_{\pm}X,Y)\theta_{\pm}+\theta_{\pm}(J_{\pm}Y)X-\theta_{\pm}(Y)J_{\pm}X$$
and therefore
\begin{equation}\label{nabla N}
\begin{array}{c}
2(\nabla_X N)Y=g(X,Y)(J_+-fJ_-)\theta_++g((J_+-fJ_-)X,Y)\theta_+\\[6pt]+ \theta_+((J_+-fJ_-)Y)X
-\theta_+(Y)(J_+-fJ_-)X +2X(f)J_-Y
\end{array}
\end{equation}
since $\theta_{-}=-\theta_{+}$. Also $p=-1/2g(J_+,J_-)=1/4tr(J_+\circ J_-)$ and we get by
(\ref{dp}) that
\begin{equation}\label{dpK}
dp=\displaystyle{\frac{1}{2}}\theta_+\circ [J_+,J_-]=\sqrt{p^2-1}\theta_+\circ K.
\end{equation}
We have $X,Y\in Ker\, N$ and
\begin{equation}\label{basis}
\begin{array}{l}
g(X,X)=g(Y,Y)=0, \\[6pt]
g(X,Y)=g(X,J_+X)=g(Y,J_+Y)=0,\\[6pt]
g(J_+X,Y)=(f^2-1)|\theta_+|^2=2(fp-1)|\theta_+|^2.
\end{array}
\end{equation}
 Then the vector fields $X,Y,J_+X, J_+Y$ form a basis of the tangent space at each point of $M$. We
have also that
\begin{equation}\label{Y}
\begin{array}{l}
g(\theta_+,
J_+J_-\theta_+)=g(\theta_+,J_-J_+\theta_+)=-g(J_+\theta_+,J_-\theta_+)=p|\theta_+|^2\\[8pt]
Y=\theta_++\displaystyle{\frac{J_+J_--2p
Id+J_+J_-}{2\sqrt{p^2-1}}\theta_+=\frac{-f\theta_++J_+J_-\theta_+}{\sqrt{p^2-1}}}.
\end{array}
\end{equation}
Now $(\nabla_X N)Y-(\nabla_Y N)X=-N\nabla_X Y+N\nabla_Y X=-N[X,Y]$ since $NX=NY=0$. To compute
$N[X,Y]$ we
 use the fact that $J_+X=-fJ_-X, J_+Y=-fJ_-Y$. Hence by (\ref{nabla N}) and (\ref{basis}) we have
 $$2|\theta_+|^{-2}(\nabla_X N)Y= 2(f^2-1)\theta_+-2J_+X=-2f\sqrt{p^2-1}Y$$ since
 $\theta_+(J_{+}Y)=g(J_{+}Y,\theta_{+})=0$ by (\ref{Y}), $\theta_+(Y)=|\theta_+|^2$ and, in view of (\ref{dpK}) and Lemma~\ref{N},  $X(f)=f\theta_{+}(KX)=f\theta_{+}(X)=fg(X,\theta_{+})=0$. Similarly
  \[
2|\theta_+|^{-2}(\nabla_Y N)X=-2(f^2-1)\theta_+ + 2(fp-1)Y-2fJ_{-}X\\
\]
\[
=2(fp+f\sqrt{p^2-1}-1)Y=0
\]
since $\theta_{+}(X)=0$,
$\theta_{+}(J_{+}X)=g(-\theta_{+}+fJ_{+}J_{-}\theta_{+},\theta_{+})=-|\theta_{+}|^2+fp|\theta_{+}|^2$
by (\ref{Y}) and $Y(f)=fg(Y,\theta_{+})=f|\theta_{+}|^2$. So $N[X,Y]=f\sqrt{p^2-1}|\theta_+|^2Y$.
We can easily check that $NJ_++J_+N=2(pf-1)Id$. Then $NJ_+Y=2(pf-1)Y$ so $[X,Y]\in Span(X,Y,J_+Y)$.
It follows from (\ref{basis}) that $X,Y,J_+Y$ are linearly independent at  every point, hence
$rank[{\cal D},{\cal D}]=3$. If $[Y,J_+Y]$ has  nowhere-vanishing $J_+X$-component, then
$rank[[{\cal D},{\cal D}],{\cal D}]=4$, so ${\cal D}$ is an Engel structure.  To find the
$J_+X$-component of $[Y,J_+Y]$ we use that $[Y,J_+Y]=\nabla_Y J_+Y-\nabla_{J_+Y}Y$. First observe
that $(\nabla_Y N)Y=0$, so $\nabla_YY\in Span\{X,Y\}$. We have also that $2(\nabla_Y J_+)Y =
-\theta_+(Y)J_+Y = -|\theta_+|^2J_+Y$. Since $\nabla_Y J_+Y = (\nabla_Y J_+)Y+J_+\nabla_Y Y$, then
$\nabla_Y J_+Y \in Span\{J_+X,J_+Y\}$. Moreover, the $J_+X$-component of $\nabla_Y J_+Y$ is equal
to the $J_+X$-component of $J_+\nabla_Y Y$ which is also the $X$-component of $\nabla_Y Y$. On the
other hand we have
\[
2(\nabla_{J_+Y} N)Y = -\theta_+(Y)(J_+-fJ_-)(J_+Y) =
2pf|\theta_+|^2 Y
\]
 since $(J_+-fJ_-)J_+Y=J_+(J_++fJ_-)Y-2pfY= -2pfY$. So $N\nabla_{J_+Y}Y=-pf|\theta_+|^2Y$ and
$\nabla_{J_+Y}Y\in Span\{X,Y,J_+Y\}$ does not have $J_+X$-component. Then $[Y,J_+Y]=\nabla_Y
J_+Y-\nabla_{J_+Y}Y$ has nowhere-vanishing $J_+X$-component iff $\nabla_Y Y$ has nowhere-vanishing
$X$-component.

 To finish the proof notice that $\nabla_YY\in Span\{X,Y\}$ and if its
 $X$-component vanishes locally, $\nabla_YY=FY$ which in turn means that the flow of $Y$ is geodesic.
{\it q.e.d.}

Note finally that if $p=const$, then $\theta_+=\theta_-$ and the distribution $\cal{D}$ is
integrable.


\vskip 20pt

\noindent Johann Davidov

\noindent Institute of Mathematics and Informatics

\noindent Bulgarian Academy of Sciences

\noindent 1113 Sofia, Bulgaria

\noindent and

\noindent "L. Karavelov" Civil Engineering Higher School

\noindent 1373 Sofia, Bulgaria

\noindent jtd@math.bas.bg,

\medskip

\noindent Oleg Mushkarov

\noindent Institute of Mathematics and Informatics

\noindent Bulgarian Academy of Sciences

\noindent 1113 Sofia, Bulgaria

\noindent muskarov@math.bas.bg

\medskip

\noindent Gueo Grantcharov and Miroslav Yotov

\noindent Department of Mathematics

\noindent Florida International University

\noindent Miami, FL 33199

\noindent grantchg@fiu.edu, yotovm@fiu.edu

\end{document}